\documentclass[11pt]{article}

\usepackage{amssymb}
\usepackage{latexsym}
\usepackage{amsmath}
\usepackage{amscd}
\usepackage{amsbsy}
\usepackage{longtable}
\usepackage[dvips]{graphics}

\usepackage[latin1]{inputenc}

\newtheorem{theorem}{Theorem}

\newtheorem{lemma}{Lemma}

\newenvironment{proof}[1][Proof]{\textbf{#1.} }{\ \rule{0.5em}{0.5em}}

\newcommand{\N}{\mathbb{N}}

\newcommand{\deq}{\mathrel{\mathop:} = } 
\begin{document}

\title{Quadratic functional estimation in inverse problems}

\author{Cristina BUTUCEA$^{1}$, Katia MEZIANI$^{2}$\\
$^{1}$ Laboratoire Paul Painlev\'e (UMR CNRS 8524), \\
Universit\'e des Sciences et Technologies de Lille 1, \\
59655 Villeneuve d'Ascq cedex, France\\
e-mail: cristina.butucea@math.univ-lille1.fr\\
$^{2}$Laboratoire de Probabilités et Modèles Aléatoires, \\
Universit\'e Paris VII (Denis Diderot),\\
75251 Paris Cedex 05, France\\
 email: meziani@math.jussieu.fr
}
\maketitle

\begin{abstract} We consider in this paper a Gaussian sequence model
of observations $Y_i$, $i\geq 1$ having mean (or signal) $\theta_i$ and variance 
$\sigma_i$ which is growing polynomially like $i^\gamma$, $\gamma >0$. 
This model describes a large panel of inverse problems.
We estimate the quadratic functional of the unknown signal $\sum_{i\geq 1}\theta_i^2$
when the signal belongs to ellipsoids of both finite smoothness functions
(polynomial weights $i^\alpha$, $\alpha>0$) and infinite smoothness 
(exponential weights $e^{\beta i^r}$, $\beta >0$, $0<r \leq 2$).
We propose a Pinsker type projection estimator in each case and study its 
quadratic risk. When the signal is sufficiently smoother than the difficulty 
of the inverse problem ($\alpha>\gamma+1/4$ or in the case of exponential weights), 
we obtain the parametric rate and the efficiency 
constant associated to it. Moreover, we give upper bounds of the second order term
in the risk and conjecture that they are asymptotically sharp minimax. 
When the signal is finitely
smooth with $\alpha \leq \gamma +1/4$, we compute non parametric upper bounds 
of the risk of and we presume also that the constant is asymptotically sharp.
\end{abstract}

\noindent {\bf Mathematics Subject Classifications 2000:} 
62F12, 62G05, 62G10, 62G20\\

\noindent {\bf Key Words: } Gaussian sequence model, 
inverse problem, minimax upper bounds, 
parametric rate, Pinsker estimator, projection estimator, 
quadratic functional, second order risk.\\


\section{Introduction}

We observe $\{Y_i\}_{i=1\cdots n}$ 
\begin{equation}
\label{model} Y_i=\theta_i+\epsilon\, \xi_i  \quad \forall i=1\cdots n
\end{equation}
where $\xi_i$ are independent identically distributed (i.i.d.) random variables, having a   
Gaussian law with zero mean and variance $\sigma_i^2 = i^{2 \gamma}$ for some fixed
$\gamma \geq 0$. Let us mention that in case $\{\sigma_i\}_{i\geq 1}$ is a bounded sequence
the problem is direct and when $\sigma_i \to \infty$ the problem is an inverse problem.
We say that the problem is ill-posed when $\sigma_i$ increases polynomially and severely 
ill-posed when it increases exponentially.

We want to estimate the quadratic functional $Q(\theta)=\sum_{i=1}^\infty \theta^2_i$, 
where $\theta =\{\theta_i\}_{i\geq 1}$ belongs to the $\ell_2$-ellipsoid
\begin{equation}
\label{class}
\Sigma=\left\{\theta: \sum_{i=1}^\infty a_i^2 \theta^2_i\leq L\right\},
\end{equation}
where $a_i$ is a non decreasing sequence of positive real numbers
and $ L>0$. We consider both polynomial sequence $a_i=i^\alpha$ 
where we say that the signal is (ordinary) smooth and exponential sequence
$a_i=\exp(\beta i^r)$ where we say that the signal is super-smooth, $\alpha,\,\beta >0$
and $0<r \leq 2$.

It is known that this model can be deduced from a linear operator
equation with noisy observations $Y=A x + \epsilon \, \xi,$ where $A:\mathcal{H} \to \mathcal{H}$
is a known linear operator on the Hilbert space $\mathcal{H}$, 
$x$ belongs to $\mathcal{H}$ is the signal of interest and 
$\xi$ is a standard white Gaussian noise.
By considering an orthonormal basis $\{\varphi_i\}_{i\geq 1}$ of $\mathcal{H}$,
we consider only the sequence of values $Y_i:= Y(\varphi_i)/b_i$, where $b_i^2$ are
the eigenvalues of $A A^*$ for $i\geq 1$.
For more details and examples of inverse problems that can be written in the form
$(\ref{model})$ we refer the reader to Cavalier~{\it et al.}$\cite{CGPT}$, $\cite{CGLT}$ and references therein.
We mention as particular examples the convolution operator, the Radon transform in the case
of tomography or problems described by partial differential equations.

\bigskip

Estimation of $\theta$ in the inverse problem $(\ref{model})$ with a quadratic risk 
was thoroughly studied in the literature from the minimax point of view. Let us only 
mention a few minimax adaptive results: oracle inequalities in
Cavalier~{\it et al.}~\cite{CGLT}, sharp adaptive estimation by block thresholding
in Cavalier~{\it et al.}~\cite{CGPT} and adaptive estimators defined by penalized 
empirical risk in Golubev~\cite{golubev} .

Estimation of quadratic functionals in inverse problems was studied in two 
particular problems (specified operators). 
Butucea~\cite{butucea} considered the convolution density model
and studied the rates of a kernel type estimator. 
Méziani~\cite{meziani} estimates the purity of a quantum state, which corresponds 
mathematically to a quadratic functional of a bivariate function of mass 1, in a 
double inverse problem: tomography and convolution with Gaussian noise.
Our model allows to consider more general inverse problems, i.e. various operators $A$.

Quadratic functionals were much more studied in the direct problem ($\sigma_j$ bounded 
for all $j$) since first results given by Ibragimov and Has'minski\u \i~\cite{ibragimov1} 
and Ibragimov {\it et al.}~\cite{ibragimov}.
Fan~\cite{fan} gave minimax rates over hyperrectangles and Sobolev-type ellipsoids. 
Donoho and Nussbaum~\cite{donohonussbaum} gave Pinsker sharp minimax estimators in 
this model and in the equivalent models of fixed equidistant design regression 
and Gaussian white noise model. Fore more general bodies which are not quadratically convex, 
Cai and Low~\cite{cailow05} showed that nonquadratic estimators attain the minimax rate
of the quadratic functional.
For adaptive estimators over hyperrectangles we cite Efromovich and Low~\cite{efrolow}. 
Sharp or nearly sharp adaptive estimators over $l_p$-bodies were found by 
Klemel\"a~\cite{klemela}. Adaptive estimators over more general Besov and $l_p$ bodies
were given by Cai and Low~\cite{cailow06}. In the density model, let us mention
adaptive estimators via model selection by Laurent~\cite{laurent}.
 
Let us underline the difference between estimating $Q(\theta)$ in our model and 
that of estimating from direct data  $\sum_{j \geq 1} j^{2\gamma} \theta_j^2$ 
for $\gamma \in \N$ as 
it was done, e.g., by Fan~\cite{fan}, Donoho and Nussbaum~$\cite{donohonussbaum}$ 
and Klemel\"a~\cite{klemela}.
In our case, the variance of our estimators is slower. When estimating the quadratic
functional of a derivative, the bias is smaller, so the rates and constants are different.

\bigskip

Here, we give a Pinsker-type projection estimator which automatically attains the
parametric rate and the efficiency constant for all super-smooth signals and
for the smooth signals when $\alpha \geq \gamma+1/4$. Moreover, in this case
we give nonparametric minimax upper bounds of the second order term in the 
quadratic risk. 
Our estimator attains the expected minimax nonparametric rate in the case 
of smooth signals with $\alpha < \gamma +1/4$. We conjecture
that the asymptotic constant in the nonparametric upper bound of the risk is sharp. 
The proofs of sharp lower bounds will make the object of future work.

Let us mention that our method can be easily adapted 
for severely ill posed inverse problems, 
i.e. $\sigma_i$ increases as an exponential. The case where $\sigma_i=e^{i^2}$ is of particular interest in practice and hasn't been studied for estimating the 
signal $\{\theta_i\}_{i\geq 1}$ either.
Future developments should concern adaptive estimation of the quadratic functional.

In Section~\ref{proc} we describe the estimator and the precise choice of tuning
parameters and give asymptotic upper bounds rates of convergence and associated constant.
We postpone the proofs to the Section~\ref{proofs} and the Appendix.


\section{Estimation procedure and results}\label{proc}

Let us define the estimator
\begin{eqnarray}
\label{estimateur}
\widetilde{Q}=\sum_{i=1}^\infty h_i (Y^2_i-\epsilon^2\sigma_i^2),
\end{eqnarray}
where $\{h_i\}_{i\geq 1}$ is a sequence between 0 and 1. We shall actually see that
the optimal sequence is truncated, i.e. $h_i=0$ for all $i > W$ and that the
optimal value of $W$ tends to infinity when $\epsilon \to  0$.

\bigskip

Let us first consider the case of smooth signal: $\theta \in \Sigma(\alpha, L)$,
where $a_i=i^{ \alpha}$.
\begin{theorem} \label{ub}
Let observations $Y_1,\ldots,Y_n,...$ satisfy model~(\ref{model}). Then
the estimator $\widetilde{Q}$ in (\ref{estimateur}) with parameters
$\{h_i\}_{i\geq 1}$ and $W$ defined by
\begin{eqnarray*}
h_i &=& \left(1 -\left(\frac i W\right)^{2\alpha} \right)_+ \text{ and }\\
W &=& \left\lfloor 
\left(\frac{L^2(4\gamma+4\alpha+1)(4\gamma+2\alpha+1)}{4\alpha} 
\right)^{\frac{1}{4\alpha+4\gamma+1}}\epsilon^{-\frac{4}{4\alpha+4\gamma+1}}
\right\rfloor
\end{eqnarray*}
is such that
\begin{eqnarray*}
\sup_{\theta \in \Sigma(\alpha,L)} E\left[\left(\widetilde{Q}- Q(\theta)\right)^2\right]
& = & C(\alpha, \, \gamma,\, L) \epsilon^{\frac{16 \alpha}{4\alpha+4\gamma+1}}
(1+o(1)), 
\end{eqnarray*}
if $\alpha \leq \gamma +\frac 14$, 
\begin{eqnarray*}
\sup_{\theta \in \Sigma(\alpha,L)} E\left[\left(\widetilde{Q}- Q(\theta)\right)^2
- 4 \epsilon^2 \sum_{i=1}^\infty\sigma_i^2\theta_i^2 \right]
& = & C(\alpha, \, \gamma, \,L) \epsilon^{\frac{16 \alpha}{4\alpha+4\gamma+1}}(1+o(1)), 
\end{eqnarray*}
if $\alpha > \gamma + \frac 14$, where
\begin{eqnarray}
C(\alpha, \, \gamma,\, L) 
&=& \frac{L^{2 \frac{4 \gamma +1}{4\alpha+ 4\gamma+1}}}{(4\gamma+1)}
\left(\frac{2\alpha + 4 \gamma +1}{4 \alpha} \right)^{-\frac{4 \alpha}{4\alpha+ 4\gamma+1}} 
(4\alpha+ 4\gamma+1)^{\frac{ 4\gamma+1}{4\alpha+ 4\gamma+1}}. \label{C1}
\end{eqnarray}
\end{theorem}

We find a known phenomenon in quadratic functional estimation literature, 
i.e. the existence of two cases: 
a regular one, where the rate is parametric $\epsilon^{-2}$, 
and an irregular case when the rate is significantly slower. 
We conjecture that Theorem~\ref{ub} exhibits sharp asymptotic constant in this last
case.

In the regular case (when the underlying signal is smoother than the 
'difficulty' of the operator $A$), Theorem~\ref{ub} says actually two things. 
One of them is that, for each $\theta$ in the set $\Sigma$ the quadratic
risk of our estimator is of parametric rate and attains the efficiency constant
in our model:
$$
E\left[\left(\widetilde{Q}- Q(\theta)\right)^2\right] 
= 4 \epsilon^2 \sum_{i=1}^\infty\sigma_i^2\theta_i^2 (1+o(1)),
$$
as $\epsilon \to 0$. Secondly, the quadratic risk is decomposed and the second order
risk is optimized for our choice of parameters and equals the risk in the 
non parametric case.

Note also, that the rates are not surprising when compared to the results of 
Butucea~\cite{butucea} for the convolution density model. 
No second order terms were evaluated there, nor constants associated to the
nonparametric rate. The efficiency constant is naturally different for the density model.

\bigskip

Let us now consider the case of super-smooth signal: $\theta \in \Sigma(\beta,\,r,\, L)$,
where $a_i=\exp( \beta i^r)$.
\begin{theorem} \label{ub2}
Let observations $Y_1,\ldots,Y_n,...$ satisfy model~(\ref{model}). Let
the estimator $\widetilde{Q}$ in (\ref{estimateur}) be defined with parameters
$\{h_i\}_{i\geq 1}$ given by
\begin{eqnarray*}
h_i &=& \left(1-\frac{e^{2\beta i^r}}{e^{2\beta W^r}}\right)_+  
\end{eqnarray*}
and $W$ solution of the equation
\begin{eqnarray*}
W^{4\gamma+(1-r)_+}\exp(4\beta W^r- 2\beta r W^{r-1}I_{(r>1)})
= c(\beta,r,\gamma,L) \epsilon^{-4},
\end{eqnarray*}
with the constant $c:=c(\beta,r,\gamma,L)=2 \beta r L^2$ if $0<r<1$, 
$c=L^2(e^{4 \beta}-1)/(2 e^{2 \beta})$ if $r=1$, 
$c=L^2 /2 $ if $1<r<2$ and $c= L^2/(2 e^{2 \beta})$ if $r=2$.
Then 
\begin{eqnarray*}
\sup_{\theta \in \Sigma(\beta,r,L)} E\left[\left(\widetilde{Q}- Q(\theta)\right)^2
- 4 \epsilon^2 \sum_{i=1}^\infty\sigma_i^2\theta_i^2 \right]
& = & \frac{2 \epsilon^4}{4\gamma +1} 
\left( \frac{\log(1/\epsilon)}{\beta}\right)^{(4 \gamma + 1)/r}(1+o(1)). 
\end{eqnarray*}

\end{theorem}

We note that in this case, the signal is always smoother than the difficulty of 
the inverse problem, so there is always a parametric rate term in the quadratic risk.
Our estimator also optimizes the upper bounds for the second order term in the quadratic
risk. In this last term, the bias term is always smaller than the variance term
for super-smooth signals.

\section{Proofs}\label{proofs}

\begin{proof}[Proof of Theorem~\ref{ub}]
We decompose as usually the quadratic risk 
$E\left[\left(\widetilde{Q}- Q(\theta)\right)^2\right]$ into bias plus variance.
The bias term can be written
\begin{eqnarray}
\left( E[\widetilde{Q}]-Q(\theta) \right)^2&=&\left(\sum_{i=1}^\infty h_i E[Y^2_i-\epsilon^2\sigma_i^2]-\sum_{i=1}^\infty  \theta_i^2 \right)^2\nonumber\\
\label{bias}
&=&\left(\sum_{i=1}^\infty h_i \theta_i^2-\sum_{i=1}^\infty  \theta_i^2 \right)^2=\left(\sum_{i=1}^\infty  \theta_i^2(1-h_i)\right)^2.
\end{eqnarray}
The variance term is decomposed as follows
\begin{eqnarray}
E\left[\left(\widetilde{Q}-E[\widetilde{Q}]\right)^2\right]
&=& E\left[\left(\sum_{i=1}^\infty h_i (Y^2_i-\epsilon^2\sigma_i^2)-\sum_{i=1}^\infty h_i\theta_i^2 \right)^2\right]\nonumber\\
&=& E\left[\left(\sum_{i=1}^\infty h_i (Y^2_i-\epsilon^2\sigma_i^2-\theta_i^2) \right)^2\right]\nonumber.
\end{eqnarray}
Since $Y_i$ are independent and $\xi_i$ are independent Gaussian random variables:
\begin{eqnarray}
E\left[\left(\widetilde{Q}-E[\widetilde{Q}]\right)^2\right]
&=&\sum_{i=1}^\infty h_i^2 E\left[(Y^2_i-\epsilon^2\sigma_i^2-\theta_i^2)^2\right]\\
&=&\sum_{i=1}^\infty h_i^2 
E\left[(2\epsilon\theta_i\xi_i-\epsilon^2\sigma_i^2+\epsilon^2\xi_i^2)^2\right]\nonumber\\
&=&\sum_{i=1}^\infty h_i^2 \left\{\epsilon^4E\left[\xi_i^4\right]-2\epsilon^4\sigma_i^2
E\left[\xi_i^2\right] +4\epsilon^2\theta_i^2E\left[\xi_i^2\right]
+\epsilon^4\sigma_i^4 \right\}\nonumber.
\end{eqnarray}

Now, use the facts that $E[\xi_i^2]=\sigma_i^2$ and $E[\xi_i^4]=3\sigma_i^4$ to get
\begin{eqnarray}
E\left[\left(\widetilde{Q}-E[\widetilde{Q}]\right)^2\right]
\label{varianceun}
&=&4\epsilon^2\sum_{i=1}^\infty h_i^2\sigma_i^2\theta_i^2+2\epsilon^4\sum_{i=1}^\infty h_i^2\sigma_i^4\nonumber\\
\label{variance}
&=&4\epsilon^2\sum_{i=1}^\infty\sigma_i^2\theta_i^2-4\epsilon^2\sum_{i=1}^\infty (1-h_i^2)\sigma_i^2\theta_i^2+2\epsilon^4\sum_{i=1}^\infty h_i^2\sigma_i^4
\end{eqnarray}
Thus by \eqref{bias} and \eqref{variance} we get
\begin{eqnarray}
\label{termerisque}
E\left[\left(\widetilde{Q}-Q(\theta)\right)^2\right]
=A_0(h,\theta)+A_1(h)+A_2(\theta)-A_3(h,\theta),
\end{eqnarray}
where 
\begin{eqnarray*}
A_0(h,\theta)=A_0&:=&\left(\sum_{i=1}^\infty\theta_i^2(1-h_i)\right)^2,\\
A_1(h)=A_1&:=&2\epsilon^4\sum_{i=1}^\infty h_i^2\sigma_i^4,\\
A_2(\theta)=A_2&:=&4\epsilon^2\sum_{i=1}^\infty\sigma_i^2\theta_i^2, \\
A_3(h,\theta) = A_3&:=&4\epsilon^2\sum_{i=1}^\infty (1-h_i^2)\sigma_i^2\theta_i^2.
\end{eqnarray*}
 If we note $T(h,\theta):=A_0(h,\theta)+A_1(h)=A_0+A_1$, then we want to find 
$$\inf_{h}\sup_{\theta\in\Sigma}T(h,\theta)
\leq\sup_{\theta\in\Sigma}T(h,\theta)
\leq\sup_{\theta\in\partial\Sigma}T(h,\theta)$$
where the infimum is taken with respect to all sequences $h$ such that $0 \leq h_i \leq 1$
for all $i\geq 1$ and with 
\begin{equation}
\label{class1}
\partial\Sigma=\left\{\theta: \sum_{i=1}^\infty a_i^2\theta^2_i= L\right\}.
\end{equation}

Let us define 
$F(h,\theta)=T(h,\theta)-\kappa\left(\sum_{i=1}^\infty a_i^2\theta^2_i-L\right)$ 
with $\kappa>0$.
Then for all $j\in\mathbb{N}^*$ the optimal $h$ and $\theta$ have to verify
\begin{eqnarray*}
\frac{\partial}{\partial\theta_j}F(h,\theta)=0\quad\text{and}\quad
\frac{\partial}{\partial h_j}F(h,\theta)=0.
\end{eqnarray*}
We get 
\begin{eqnarray}
\label{optimisation}
h_j&=&\left(1-\frac{\kappa a_j^2}{2\sum_{i=1}^\infty\theta_i^2(1-h_i)}\right)_+
=\left(1-\widetilde{\kappa} a_j^2\right)_+,\nonumber\\
(\theta_j^*)^2&=&\frac{2\epsilon^4\sigma_j^4 h_j}{\sum_{i=1}^\infty\theta_i^2(1-h_i)},
\end{eqnarray}
where $\widetilde{\kappa}>0$.
Let us write $h_j= \left(1-\frac{j^{2\alpha}}{W^{2\alpha}}\right)_+$ where $W \to \infty$
when $\epsilon \to 0$.

Recall that $\Sigma=\Sigma(\alpha,L)=\left\{\theta: \sum_{i=1}^\infty i^{2\alpha} \theta^2_i\leq L\right\}$ then for $\theta^*\in\partial\Sigma(\alpha,L)$ 
we can write both
\begin{eqnarray*}
\sum_{i=1}^\infty\theta_i^{*2}(1-h_i)
&=&\frac{1}{W^{2\alpha}}\sum_{i=1}^W\theta_i^{*2}i^{2\alpha}+\sum_{i>W}\theta_i^{*2}\\
&\leq &\frac{1}{W^{2\alpha}}\sum_{i=1}^W\theta_i^{*2}i^{2\alpha}
+\frac{1}{W^{2\alpha}}\sum_{i>W}\theta_i^{*2}i^{2\alpha}\leq\frac{L}{W^{2\alpha}}
\end{eqnarray*}
and
\begin{eqnarray*}
\sum_{i=1}^\infty\theta_i^{*2}(1-h_i)&\geq&\frac{1}{W^{2\alpha}}\sum_{i=1}^W\theta_i^{*2}i^{2\alpha}\\
&\geq&\frac{1}{W^{2\alpha}}\sum_{i=1}^\infty\theta_i^{*2}i^{2\alpha}-\frac{1}{W^{2\alpha}}\sum_{i>W}\theta_i^{*2}i^{2\alpha}
=\frac{L}{W^{2\alpha}}(1-o(1)).
\end{eqnarray*}
Therefore $A_0 = L^2 W^{-4 \alpha} (1+o(1))$, as $\epsilon \to 0$. This means also
that we can write
$$
(\theta_j^*)^2=\frac{2\epsilon^4\sigma_j^4 W^{2 \alpha}}{L} 
\left(1-\frac{j^{2 \alpha}}{W^{2 \alpha}}\right)_+.
$$

Let us now compute the optimal $W$, using again the fact that 
$\theta^*\in\partial\Sigma(\alpha,L)$ which is equivalent to
\begin{eqnarray*}
&&\sum_{i=1}^\infty i^{2\alpha}(\theta_i^*)^2= L.
\end{eqnarray*}
This is further equivalent to
\begin{eqnarray*}
W^{4\alpha+4\gamma+1}\frac{1}{W}\sum_{i=1}^{W} 
\left(\frac{i}{W}\right)^{4\gamma+2\alpha}\left(1-\left(\frac{i}{W}\right)^{2\alpha} \right)
=\frac{L^2}{2\epsilon^4}
\end{eqnarray*}
giving
\begin{eqnarray*}
&& \frac{2\alpha W^{4\alpha+4\gamma+1}}{(4\gamma+4\alpha+1)(4\gamma+2\alpha+1)}(1+o(1))
=\frac{L^2}{2\epsilon^4}.
\end{eqnarray*}
Therefore
\begin{eqnarray}
\label{M}
W=\left(\frac{L^2}{B(\alpha,\gamma)} \right)^{\frac{1}{4\alpha+4\gamma+1}}\epsilon^{-\frac{4}{4\alpha+4\gamma+1}}(1+o(1)),
\end{eqnarray}
where $B(\alpha,\gamma)\deq \frac{4\alpha}{(4\gamma+4\alpha+1)(4\gamma+2\alpha+1)}$ 
and we'll take $W$ to be the integer part of the dominant term. From now on, 
we denote $B\deq B(\alpha,\gamma)$.

We have to evaluate the term  defined in \eqref{termerisque}. 
For $\alpha\leq \gamma+\frac 14$, we have 
\begin{eqnarray*}
\label{A0}
A_0&=&\left(\sum_{i=1}^\infty\theta_i^2(1-h_i)\right)^2=L^2 W^{-4\alpha}(1+o(1))\nonumber\\
&=&\left(L^{2(4\gamma+1)}B^{4\alpha}\right)^{\frac{1}{4\gamma+4\alpha+1}}\epsilon^{\frac{16\alpha}{4\gamma+4\alpha+1}}(1+o(1)),\\
\label{A1}
A_1&=&2\epsilon^4\sum_{i=1}^\infty \sigma_i^4h_i^{2}=2\epsilon^4 W^{4\gamma+1}\frac{1}{W}\sum_{i=1}^{W}\left(\frac{i}{W}\right)^{4\gamma}\left(1-\left(\frac{i}{W}\right)^{2\alpha}\right)^2\nonumber\\
&=&\frac{16\alpha^2\epsilon^4W^{4\gamma+1}}{(4\gamma+1)(4\gamma+4\alpha+1)(4\gamma+2\alpha+1)}(1+o(1))\nonumber\\
&=&\frac{4\alpha}{4\gamma+1}\left(L^{2(4\gamma+1)}B^{4\alpha}\right)^{\frac{1}{4\gamma+4\alpha+1}}\epsilon^{\frac{16\alpha}{4\gamma+4\alpha+1}}(1+o(1)),\\
\label{A2}
A_2&=&4\epsilon^2\sum_{i=1}^\infty\sigma_i^2\theta_i^{*2}=\frac{8\epsilon^6 W^{6\gamma+2\alpha+1}}{L} \frac{1}{W}\sum_{i=1}^{W}\left(\frac{i}{W}\right)^{6\gamma}\left(1-\left(\frac{i}{W}\right)^{2\alpha}\right)\nonumber\\
&=&\frac{16\alpha\epsilon^6W^{6\gamma+2\alpha+1}}{L(6\gamma+1)(6\gamma+2\alpha+1)}(1+o(1))\nonumber\\
&=&O(1)\epsilon^{\frac{16\alpha+2}{4\alpha+4\gamma+1}}(1+o(1))=o(1)A_1,
\end{eqnarray*}
as $\epsilon \to 0$.
As $h_i\in[0,1]$ for all $i\in\mathbb{N}$, the term $A_3=4\epsilon^2\sum_{i=1}^\infty (1-h_i^2)\sigma_i^2\theta_i^2\leq A_2$. Then the quadratic risk is such that
\begin{eqnarray*}
E\left[\left(\widetilde{Q}- Q(\theta)\right)^2\right]&=&\left(A_0+A_1\right)(1+o(1))\\
&=&\left(L^{2(4\gamma+1)}B^{4\alpha}\right)^{\frac{1}{4\gamma+4\alpha+1}}\frac{4\gamma+4\alpha+1}{4\gamma+1}\epsilon^{\frac{16\alpha}{4\gamma+4\alpha+1}}(1+o(1)),
\end{eqnarray*}
as $\epsilon \to 0$ and this explains the constant $C(\alpha, \gamma, L)$ in (\ref{C1}).

Let us note that if  $\alpha>\gamma+\frac 14$, we can estimate the quadratic 
functional at the parametric rate as $A_2$ is the dominant term in the risk 
and is of order $\epsilon^2$. More precisely
\begin{eqnarray*}
E\left[\left(\widetilde{Q}- Q(\theta)\right)^2\right]
&=& 4 \epsilon^2 \sum_{i=1}^\infty \sigma_i^2 \theta_i^2 \, (1+o(1)) = A_2 (1+o(1)),
\end{eqnarray*}
as $\epsilon \to 0$.
Indeed, it is easy to see that in this case 
$$
A_0+A_1 = C(\alpha, \gamma, L) \epsilon^{\frac{16\alpha}{4\gamma+4\alpha+1}}(1+o(1))
= o(A_2)
$$
and, moreover,
\begin{eqnarray*}
A_3
&=& 4\epsilon^2 \sum_{i=1}^{W}\left[1-\left(1-\frac{i^{2 \alpha}}{W^{2\alpha}}\right)^2\right] i^{2\gamma} \theta_i^2+ 4 \epsilon^2 \sum_{i> W}  i^{2\gamma} \theta_i^2\\
& \leq & 4 \epsilon^2 \sum_{i=1}^{W}\frac{i^{2 \alpha+2\gamma}}{W^{2\alpha}}+ 4\epsilon^2 \sum_{i> W} i^{2(\gamma -\alpha)} i^{2\alpha} \theta_i^2\\
&\leq & 4 \epsilon^2 W^{2(\gamma-\alpha)}\sum_{i=1}^{W}\left(\frac{i}{W}\right)^{2\gamma}i^{2\alpha} \theta_i^2+4\epsilon^2 W^{2(\gamma-\alpha)}\sum_{i=1}^{W}i^{2\alpha} \theta_i^2\\
&\leq & 4 \epsilon^2 W^{2(\gamma-\alpha)} L 
= O(1) \epsilon^{\frac{16\alpha+2}{4\alpha+4\gamma+1}} = o(A_0+A_1), 
\end{eqnarray*}
as $\epsilon \to 0$.
\end{proof}

\begin{proof}[Proof of Theorem~\ref{ub2}] 
We follow the lines of proof of  Theorem~\ref{ub}. In this case, there is always
a parametric term and we do the computations of the second order term in the
quadratic risk.

We solve the same optimisation problem and find
\begin{eqnarray}
\label{parameters1}
h_i =\left(1-\frac{e^{2\beta i^r}}{e^{2\beta W^r}}\right)_+ \quad 
(\theta_j^*)^2&=&\frac{2\epsilon^4\sigma_j^4 h_j}{\sum_{i=1}^\infty\theta_i^2(1-h_i)}.
\end{eqnarray}
Then for $\theta^*\in\partial\Sigma(\beta,L,r)$ we get
\begin{eqnarray*}
\sum_{i=1}^\infty\theta_i^{*2}(1-h_i)&=&\frac{1}{e^{2\beta W^r}}\sum_{i=1}^W e^{2\beta i^r}\theta_i^{*2}+\sum_{i>W} \theta_i^{*2}\\
&\leq&\frac{1}{e^{2\beta W^r}}\sum_{i=1}^W e^{2\beta i^r}\theta_i^{*2}+\frac{1}{e^{2\beta W^r}}\sum_{i>W}e^{2\beta i^r}\theta_i^{*2}
=\frac{L}{e^{2\beta W^r}}
\end{eqnarray*}
and
\begin{eqnarray*}
\sum_{i=1}^\infty\theta_i^{*2}(1-h_i)
&\geq&\frac{1}{e^{2\beta W^r}}\sum_{i=1}^W e^{2\beta i^r}\theta_i^{*2}\\
&=&\frac{1}{e^{2\beta W^r}}\sum_{i=1}^\infty e^{2\beta i^r}\theta_i^{*2}-\frac{1}{e^{2\beta W^r}}\sum_{i>W} e^{2\beta i^r}\theta_i^{*2}
=\frac{L}{e^{2\beta W^r}}(1-o(1)).
\end{eqnarray*}
Therefore
\begin{eqnarray*}
A_0=L^2 e^{-4\beta W^r}(1+o(1)),\text{ as }\epsilon\rightarrow 0.
\end{eqnarray*}
By \eqref{parameters1}, this gives $\theta_i^{*2}=\frac{2\epsilon^4\sigma_j^4}{L}\left(e^{2\beta W^r} -e^{2\beta j^r} \right)_+ .$

To compute optimal $W$, we also use the fact $\theta^*\in\partial\Sigma(\beta,L,r)$.

\begin{eqnarray*}
&&\sum_{i=1}^\infty e^{2\beta i^r}(\theta_i^*)^2= L\Leftrightarrow\,e^{2\beta W^r}\sum_{i=1}^{W-1} i^{4\gamma}e^{2\beta i^r}-\sum_{i=1}^{W-1} i^{4\gamma}e^{4\beta i^r}=\frac{L^2}{2\epsilon^4}
\end{eqnarray*}
By using Lemmata~\ref{lm1bis} and~\ref{somme}, we have $W$ solution of the following equation
\begin{equation}
\label{Parameters3}
            \begin{array}{ll}
             W^{4\gamma} e^{ 4 \beta W^r - 2 \beta r W^{r-1}}
=c \epsilon^{-4},& \text{ if } 1<r \leq 2,\\
             W^{4\gamma} e^{4\beta W}= c \epsilon^{-4},&\text{ if } r=1,\\
             W^{4\gamma-r+1}e^{4\beta W^r}= c \epsilon^{-4},& \text{ if } 0<r<1,
        \end{array}
    \end{equation}
as $\epsilon\rightarrow 0$, with the constant $c=c(\beta, \gamma, L)$ 
defined in Theorem~\ref{ub2}. 

We evaluate $A_0+A_1$: in each of the previous cases, 
the bias term $A_0$ is infinitely smaller
than the variance term $A_1$ and the main term in $A_1$ can be given for 
$$
W=\left( \frac{\log(1/\epsilon)}{\beta}\right)^{1/r}.
$$
Indeed, by using Lemmata~\ref{lm1bis} and~\ref{somme},
\begin{eqnarray*}
A_1&=& 2\epsilon^4\sum_{i=1}^\infty \sigma_i^4h_i^{2}
=2\epsilon^4\sum_{i=1}^{W}i^{4\gamma}
\left(1-\frac{e^{2\beta i^r}}{e^{2\beta W^r}}\right)^2\nonumber\\
&=&\frac{2\epsilon^4W^{4\gamma+1}}{4\gamma+1}(1+o(1))
=\frac{2\epsilon^4 }{4\gamma+1}
\left(\frac{\log(1/\epsilon)}{\beta}\right)^{(4\gamma+1)/r}(1+o(1))=o(A_2).\\
\end{eqnarray*}
As $A_0= o(A_1)$ it is easy to see that in this case 
$$
A_0+A_1 = \frac{2\epsilon^4 }{4\gamma+1}
\left(\frac{\log(1/\epsilon)}{\beta}\right)^{(4\gamma+1)/r}(1+o(1))
= o(A_2)
$$
as $\epsilon \to 0$.

The last thing to check is that $A_3=o(A_0+A_1)$ as $\epsilon \to 0$:
\begin{eqnarray*}
A_3&=&4\epsilon^2\sum_{i=1}^\infty(1-h_i^{2})\sigma_i^2\theta_i^{*2}
\leq  8\epsilon^2 \sum_{i=1}^{W}\frac{e^{2\beta i^r}}{e^{2\beta W^r}}i^{2\gamma}\theta_i^{*2}+ 4\epsilon^2 \sum_{i> W} i^{2\gamma}\theta_i^{*2}\nonumber\\
&\leq & 8\epsilon^2 \frac{W^{2\gamma}}{e^{2\beta W^r}}\sum_{i=1}^{W}e^{2\beta i^r}\theta_i^{*2}+ 8\epsilon^2 \sum_{i> W} i^{2\gamma}\frac{e^{2\beta i^r}}{e^{2\beta i^r}}\theta_i^{*2}\nonumber\\
&\leq & 8\epsilon^2 \frac{W^{2\gamma}}{e^{2\beta W^r}}\sum_{i=1}^{W}e^{2\beta i^r}\theta_i^{*2}+ 8\epsilon^2 \frac{W^{2\gamma}}{e^{2\beta W^r}}\sum_{i> W} e^{2\beta i^r}\theta_i^{*2}\nonumber\\
&= & 8\epsilon^2 \frac{W^{2\gamma}}{e^{2\beta W^r}}L
=O(1) W^{4 \gamma +1}\epsilon^4 \frac{1}{W^{2\gamma+1 }\epsilon^2e^{2\beta W^r}}.
\end{eqnarray*}
So, we can write that
$$
A_3 = O(A_1) \frac{1}{W^{2\gamma+1 }\epsilon^2e^{2\beta W^r}}.
$$
By \eqref{Parameters3}, we easily see that
\begin{equation*}
            \begin{array}{ll}
             W^{2\gamma+1} e^{2 \beta W^r -\beta r W^{r-1}}\epsilon^{2}
=\sqrt{c}\, W,& \text{ if } 1< r \leq 2,\\
             W^{2\gamma+1} e^{2\beta W}\epsilon^{2}= \sqrt{c}\, W,& \text{ if } r=1,\\
             W^{2\gamma+1} e^{2\beta W}\epsilon^{2}= \sqrt{c} \, W^{(1+r)/2},
& \text{ if } 0<r<1,
        \end{array}
    \end{equation*}
    Then, as $W \to \infty$, we get 
for all $r\in]0,2]$, $A_3=o(A_1)$ as $\epsilon\rightarrow 0$.

\end{proof}

\section{Appendix}

\begin{lemma}
\label{lm1bis}
For all $a,\, b,\,s\,>0$ and $ v>0$
\begin{eqnarray*}
\int_0^{v} x^{a}e^{b x^s}dx&=&\frac{v^{a-s+1}e^{b v^s}}{bs}(1+o(1)),
\end{eqnarray*}as $v\rightarrow\infty$.
\end{lemma}

\begin{lemma}
\label{somme}
For $a\geq 0$, $b>0$, and $r>0$ as $N \to \infty$
\begin{equation*}
\sum_{i=1}^N i^a e^{b i^r}=\left\lbrace
  \begin{array}{c l}
      & N^a e^{b N^r}(1+o(1))\text{ if }r>1,\\
         &  \frac 1{br}N^{a+1-r}e^{bN^r}(1+o(1))\text{ if }0<r<1,\\
         & \frac{1}{(e^{b}-1)}N^a e^{b (N+1)}(1+o(1))\text{ if r=1 and }a\neq 0 .
  \end{array}
\right. 
\end{equation*} 
\end{lemma}
\begin{proof}[\textit{\textbf{Proof of Lemma~\ref{somme}}}]
$\bullet$ When $ r>1$
\begin{eqnarray*}
\sum_{i=1}^N i^a e^{b i^r}- N^a e^{b N^r} 
&=& \sum_{i=1}^{N-1} i^a e^{b i^r}\leq (N-1)^{a+1} e^{b (N-1)^r}\\
&\leq & N^a e^{bN^r} O(N)e^{-brN^{r-1}}=o(1)N^a e^{bN^r},
\end{eqnarray*}
as $N \to \infty$.

$\bullet$ When $0< r<1$
\begin{eqnarray*}
\int_1^{N+1} x^a e^{b x^r}dx\geq \sum_{i=1}^N i^a e^{b i^r}\geq \int_0^N x^a e^{b x^r}dx.
\end{eqnarray*}
Use Lemma~\ref{lm1bis} and the fact that
\begin{eqnarray*}
\int_1^{N+1} x^a e^{b x^r}dx = \int_0^N x^a e^{b x^r}dx (1+o(1).
\end{eqnarray*}

$\bullet$ When $ r=1$ we write both
\begin{eqnarray*}
\sum_{i=1}^N i^a e^{b i}&=& N^a e^{bN}+\sum_{i=1}^{N-1}i^a e^{b i}
\end{eqnarray*}
and
\begin{eqnarray*}
\sum_{i=1}^N i^a e^{b i}&=& e^b \sum_{i=0}^{N-1}(i+1)^a e^{b i}
=e^b + e^b \sum_{i=1}^{N-1}(i+1)^a e^{b i}.
\end{eqnarray*}
As the sums $\sum_{i=1}^{N-1}i^a e^{b i}$ and $\sum_{i=1}^{N-1}(i+1)^a e^{b i}$
have equivalent general terms and diverge, than they are equivalent to $S_{N-1}$, say.
We get that, for large $N$,
$$
S_N = \frac{N^a e^{b(N+1)}}{e^b -1}(1+o(1)).
$$
\end{proof}

\bibliographystyle{acm}

\bibliography{bm}

\end{document}